\documentclass[12pt]{amsart}
\usepackage[colorlinks=true,citecolor=magenta,linkcolor=magenta]{hyperref}
\usepackage{amssymb,verbatim,amscd,amsmath,graphicx}
\usepackage{enumerate}
\usepackage{graphicx}
\usepackage{caption}
\usepackage{subcaption}
\usepackage{pdfsync}
\usepackage{color,colortbl}
\usepackage{fullpage}
\usepackage{bbm}
\usepackage{comment}

\numberwithin{equation}{section}
\newtheorem{theorem}{Theorem}[section]
\newtheorem{lemma}[theorem]{Lemma}
\newtheorem{proposition}[theorem]{Proposition}
\newtheorem{corollary}[theorem]{Corollary}

\theoremstyle{definition}

\newtheorem{conjecture}[theorem]{Conjecture}
\newtheorem{def-prop}[theorem]{Definition-Proposition}

\newtheorem{remark}[theorem]{Remark}
\newtheorem{example}[theorem]{Example}

\newtheorem*{acknowledgement}{Acknowledgements}

\newtheorem{question}[theorem]{Question}

\DeclareMathOperator{\depth}{depth}

\DeclareMathOperator{\Ass}{Ass}
\DeclareMathOperator{\Min}{Min}

\DeclareMathOperator{\pd}{pd}

\newcommand{\ZZ}{{\mathbb Z}}
\newcommand{\NN}{{\mathbb N}}

\newcommand{\kk}{{\mathbbm k}}

\def\mm{{\frak m}}

\def\pp{{\frak p}}
\def\qq{{\frak q}}

\def\a{{\bf a}}

\def\1{{\bf 1}}
\def\0{{\bf 0}}

\begin{document}
	
\title{Depth functions and symbolic depth functions\\ of homogeneous ideals}

\author{Huy T\`ai H\`a}
\address{Tulane University \\ Department of Mathematics \\
	6823 St. Charles Ave. \\ New Orleans, LA 70118, USA}
\email{tha@tulane.edu}

\author[Trung]{Ngo Viet Trung}
\address{International Centre for Research and Postgraduate Training\\ 
	Institute of Mathematics \\ 
	Vietnam Academy of Science and Technology \\ 
	18 Hoang Quoc Viet, Hanoi, Vietnam}
\email{nvtrung@math.ac.vn}

\thanks
{The first author acknowledges supports from Louisiana Board of Regents, grant \# LEQSF(2017-19)-ENH-TR-25. The second author is partially supported by grant 101.04-2019.313 of the Vietnam National Foundation for Science and Technology Development.}


\begin{abstract}
We survey recent studies and results on the following problem: for which function $f: \NN \rightarrow \ZZ_{\ge 0}$ does there exist a homogeneous ideal $Q$ in a polynomial ring $S$ such that (a) $\depth S/Q^t = f(t)$ for all $t \ge 1$, or (b) $\depth S/Q^{(t)} = f(t)$ for all $t \ge 1$?
\end{abstract}

\maketitle

\centerline{\em \small Dedicated to David Eisenbud on the occasion of his 75th birthday}


\section{Introduction} \label{sec.intro}

Let $\kk$ be a field and let $S$ be a standard graded $\kk$-algebra. For a homogeneous ideal $Q \subseteq S$, we call the functions $\depth S/Q^t$ and $\depth S/Q^{(t)}$, for $t \ge 1$, the \emph{depth function} and the \emph{symbolic depth function} of $Q$, respectively.  

Depth is an important cohomological invariant (cf. \cite{AB1956, BH1993, Serre1965}). For instance, we can compute the projective dimension via depth by the Auslander and Buchsbaum formula:
$$\pd S/Q = \dim S - \depth S/Q.$$
However, our understanding of the depth function and the symbolic depth function of ideals has been quite limited. This is partly because there are no effective methods to compute and/or to compare the depth of powers and symbolic powers of an arbitrary ideal. The aim of paper is to present recent studies, which have led to satisfactory solutions to the problem of classifying depth functions and symbolic depth functions of homogeneous ideals in polynomial rings.

It is a classical result of Brodmann \cite{Brod1979} that the depth function of an ideal in a Noetherian ring is asymptotically a constant function. From a few initially known examples, the depth function of an ideal appeared to be a non-increasing function. As more examples surfaced, it became a surprising fact that the depth function may otherwise exhibit wild behaviors; see \cite{BHH2014, HTT2016,HH2005,MST2018}.  Herzog and Hibi \cite{HH2005} conjectured that the eventual behavior, as shown by Brodmann's result, is the only condition for the depth function of homogeneous ideals in polynomial rings.
In Section \ref{sec.depth}, we survey results in our recent joint work with H.D. Nguyen and T.N. Trung \cite{HNTT2020PAMS}, in which we establish the conjecture of Herzog and Hibi in its full generality.

Symbolic depth functions are much less understood compared to depth functions. For instance, there is no similar result to that of Brodmann for the eventual behavior of symbolic depth functions of ideals. This is because the symbolic Rees algebra $R_s(Q) := \bigoplus_{t \ge 0}Q^{(t)}$, which governs the behavior of symbolic powers of $Q$, is not always finitely generated. If $R_s(Q)$ is finitely generated then $\depth S/Q^{(t)}$ is an asymptotically periodic function. 
In Section \ref{sec.symbolic}, we survey recent results of the second author and H.D. Nguyen in \cite{HopT2019a, HopT2019b}, which shows that any positive and asymptotically periodic numerical function is the symbolic depth function of a homogeneous ideal in a polynomial ring.

In both Sections \ref{sec.depth} and \ref{sec.symbolic}, we shall thoroughly explain the ideas and techniques which have led to the surveyed results in \cite{HNTT2020PAMS} and \cite{HopT2019a, HopT2019b}. We believe that {\em they may provide effective tools for the study of other numerical invariants}, such as the projective dimension, the Castelnuovo-Mumford regularity and the number of associated primes of powers and symbolic powers of homogeneous ideals.

We end the paper with Section \ref{sec.question}, where we discuss a number of open questions on depth functions and symbolic depth functions, and related problems on the projective dimension of powers and symbolic powers of homogeneous ideals. For unexplained notions and terminology we refer the readers to \cite{BH1993}.

\begin{acknowledgement}
	We thank Irena Peeva for inviting us to contribute a paper to this special volume.
\end{acknowledgement}


\section{Ordinary depth functions} \label{sec.depth}

One of the main motivations for the study of depth functions of ideals is the following classical result of Brodmann \cite{Brod1979}.

\begin{theorem} \label{thm.Brod} \cite[Theorem (2)]{Brod1979}	
	Let $S$ be a Noetherian ring and let $Q \subseteq S$ be an ideal. Then, $\depth S/Q^t$ is asymptotically a constant function, i.e., $\depth S/Q^t = \depth S/Q^{t+1}$ for all $t \gg 0$.
\end{theorem}

The first systematic study on depth functions of homogeneous ideals was carried out by Herzog and Hibi \cite{HH2005}. In their work, Herzog and Hibi observed that Theorem \ref{thm.Brod} is only a special case of a more general phenomenon, which we shall now describe.
Note that a graded $S$-algebra $R$ is said to be \emph{standard graded} if it is generated over $S$ by homogeneous elements of degree one. For a graded module $E$, we denote by $E_t$ its degree $t$ component.

\begin{theorem} \label{HH1} \cite[Theorem 1.1]{HH2005}
Let $R$ be a finitely generated standard graded $S$-algebra. Let $E$ be a finitely generated graded $R$-module. Then,
$\depth E_t = \depth E_{t+1}$ for all $t \gg 0$.
\end{theorem}

For an ideal $Q \subseteq S$, let $R(Q) := \bigoplus_{t \ge 0}Q^t$ be the \emph{Rees algebra} of $Q$. Then, $R(Q)$ is a finitely generated standard graded $S$-algebra. Thus, Theorem \ref{HH1} applies to imply Theorem \ref{thm.Brod}.

Thanks to Theorem \ref{thm.Brod}, to investigate all possible depth functions, we only need to focus on \emph{convergent} non-negative \emph{numerical} functions, i.e. functions $f: \NN \rightarrow \ZZ_{\ge 0}$ with the property that $f(t) = f(n+1)$ for $t \gg 0$. 

Herzog and Hibi \cite[Theorem 4.1]{HH2005} showed that any non-decreasing convergent non-negative numerical function is the depth function of a monomial ideal. This result was astonishing because, as mentioned, depth functions of homogeneous ideals initially tend to be non-increasing. Also in \cite{HH2005},
Herzog and Hibi exhibited monomial ideals whose depth functions display unusual behaviors. These results seemed to suggest that the beginning of the depth function of a homogeneous deal in a polynomial ring can be \emph{arbitrarily} wild. In fact, at the end of \cite{HH2005}, Herzog and Hibi made the following conjecture.

\begin{conjecture}[Herzog-Hibi]
	\label{conj.HH}
	Let $f$ be any convergent non-negative numerical function. There exists a homogeneous ideal $Q$ in a polynomial ring $S$ such that $f$ is the depth function of $Q$, i.e., $f(t) = \depth S/Q^t$ for all $t \gg 1$.
\end{conjecture}

In a later work, Bandari, Herzog and Hibi  \cite{BHH2014} showed that the depth function can have an arbitrary number of local maxima and local minima, which provided a strong evidence for Conjecture \ref{conj.HH}. Until then, the constructions of particular depth functions were all more or less ad hoc. 

There was a more general attempt by the authors, together with T.N. Trung, in \cite[Theorem 4.9]{HTT2016} to show that any non-increasing convergent non-negative numerical function is the depth function of a monomial ideal. An important new ingredient in \cite{HTT2016} is the following result on the depth function of sums of ideals.

Let $A$ and $B$ be polynomial rings over a field $\kk$ with disjoint sets of variables. Let $I \subseteq A$ and $J \subseteq B$ be nonzero proper homogeneous ideals. By abuse of notations, we shall also use $I$ and $J$ to denote their extensions in the tensor product $R := A \otimes_k B$. 

\begin{proposition} \cite[Corollary 3.6(i)]{HTT2016} 
	\label{HTT2016}
Assume that $\depth I^{i-1}/I^i  \ge \depth I^i/I^{i+1}$ for $i \le t-1$. Then,
$$\depth R/(I+J)^t = \min_{i+j = n-1}\{\depth I^i /I^{i+1}+ \depth J^j/J^{i+1}\}.$$
\end{proposition}

The proof of \cite[Theorem 4.9]{HTT2016} contained an error, discovered by Matsuda, Suzuki and Tsuchiya in \cite{MST2018}. It only gave the desired conclusion for a large class of non-increasing convergent non-negative numerical functions, as established in \cite[Theorem 2.1]{MST2018}. However, a modification of this approach has finally led to a complete characterization of depth functions of monomial ideals, which confirms Conjecture \ref{conj.HH}.

\begin{theorem} \cite[Theorem 4.1]{HNTT2020PAMS} 
	\label{thm.HNTT}
	Let $f: \NN \rightarrow \ZZ_{\ge 0}$ be any convergent non-negative numerical function and let $\kk$ be any field. There exists a monomial ideal $Q$ in a polynomial ring $S$ over $\kk$ such that $f$ is the depth function of $Q$.
\end{theorem}

The key idea in \cite{HNTT2020PAMS} is that depth functions are \emph{additive}, i.e., the sum of two depth functions is again a depth function. This makes use of the following result of Hoa and Tam \cite{HT2010}. 

\begin{lemma} \cite[Lemmas 1.1 and 2.2]{HT2010}  \label{HoaTam} 
Let $A$ and $B$ be polynomial rings over $\kk$ with disjoint sets of variables. Let $I \subseteq A$ and $J \subseteq B$ be nonzero proper homogeneous ideals, which are also seen as their extensions in $R = A \otimes_\kk B$. Then, \par
{\rm (i)} $I \cap J = IJ$, and \par
{\rm (ii)} $\depth R/IJ = \depth A/I + \depth B/J + 1.$
\end{lemma}

By setting $S = R/(x-y)$ and $Q = (IJ,x-y)/(x-y)$, where $x$ and $y$ are arbitrary variables in $A$ and $B$, respectively, Lemma \ref{HoaTam} gives rise to the following result.

\begin{proposition} \cite[Proposition 2.3]{HNTT2020PAMS} \label{additivity}
Let $I \subseteq A$ and $J \subseteq B$ be homogeneous ideals as in Lemma \ref{HoaTam}. 
There exists a homogeneous ideal $Q$ in a polynomial ring $S$ such that for all $t \ge 1$,
$$\depth S/Q^t = \depth A/I^t + \depth B/J^t.$$
Moreover, if $I$ and $J$ are monomial ideals then $Q$ can be chosen to be a monomial ideal.
\end{proposition}

We will use the additivity of depth functions to construct a monomial ideal whose depth function is any given
convergent non-negative numerical function. The construction is based on the following simple arithmetic observation.

To ease on notations, we shall identify a numerical function $f: \NN \rightarrow \ZZ$ with the sequence of its values $f(1),f(2),\dots$.
Let $f$ be a convergent non-negative numerical function which is not the constant function $0,0,\dots$. 
Then $f$ can be written as a sum of numerical functions of the following two types:\par
	\hspace*{3ex}{\rm Type I:} $0,\dots,0,1,1,\dots$ \par
	\hspace*{2ex}{\rm Type II:} $0,\dots,0,1,0,0,\dots$. \par
\noindent Note that if $f$ is the constant function $0,0,\dots$
then $f$ is the depth function of the maximal homogeneous ideal in any polynomial ring over $\kk$.

By Proposition \ref{additivity} and the above observation, to prove Theorem \ref{thm.HNTT} we only need to construct monomial ideals that admit any function of both Types I and II as their depth functions. Functions of Type I are non-decreasing convergent functions, and so they are the depth functions of monomial ideals, as constructed in \cite[Theorem 4.1]{HH2005}. For functions of Type II, we shall make use of monomial ideals whose depth functions are of the form $1,\dots,1,0,0,\dots$, which exist as shown by \cite[Example 4.10]{HTT2016} and \cite[Proposition 1.5]{MST2018}. 

As before, let $A$ and $B$ be polynomial rings over $\kk$ with disjoint sets of variables, and let $R = A \otimes_\kk B$. Let $I \subseteq A$ and $J \subseteq B$ be monomial ideals with depth functions $0,\dots,0,1,1,\dots$ and $1,\dots,1,0,0,\dots$, where the first 1 of the former function and the last 1 of the later function are at the same position.  
By Proposition \ref{additivity}, the function
$\depth R/((IJ)^t,x-y)$ is of the form $1,\dots,1,2,1,1,\dots$ for some variables $x,y$.
If we can find variables $x'$ and $y'$ such that $x'-y'$ is a non-zerodivisor in $R/((IJ)^t,x-y)$ for all $t \ge1$,
then 
$$\depth R/((IJ)^t,x-y,x'-y') = \depth R/((IJ)^t,x-y)-1$$ 
is of the form $0,\dots,0,1,0,0,\dots$, i.e., of Type II.
Clearly, we can identify $S = R/(x-y,x'-y')$ with a polynomial ring
and $(IJ,x-y,x'-y')/(x-y,x'-y')$ with a monomial ideal in $S$.
To find such variables $x'$ and $y'$ we need to know the associated primes of the ideal $((IJ)^t,x-y)$ for all $t \ge 1$. This is given in the next proposition.

For an ideal $Q$, denote the set of the associated primes and the set of the minimal associated primes of $Q$ by $\Ass(Q)$ and $\Min(Q)$, respectively.

\begin{proposition} \cite[Proposition 3.2]{HNTT2020PAMS} \label{control}
Let $I \subseteq A$ and $J \subseteq B$ be proper monomial ideals in polynomial rings. Let $x$ and $y$ be variables in $A$ and $B$, respectively. Then, $\Ass(IJ,x-y)$ is given by
\begin{align*}
\{(\pp,x-y)|\ \pp \in \Ass (I)\} \cup \{(\qq,x-y)|\ \qq \in \Ass(J)\} \cup \left(\bigcup_{\begin{subarray}{l} \pp \in \Ass (I), x \in \pp\\ \qq \in \Ass (J), y \in \qq \end{subarray}} \Min(\pp+\qq)\right).
\end{align*}
\end{proposition} 

Using Proposition \ref{control} one can give sufficient conditions for the existence of variables $x',y'$ such that $x'-y'$ is a non-zerodivisor in $R/((IJ)^t,x-y)$ for all $t \ge 1$ \cite[Proposition 3.5]{HNTT2020PAMS}.  
It turns out that the monomial ideals $I$ and $J$, as exhibited in \cite[Theorem 4.1]{HH2005} and \cite[Example 4.10]{HTT2016}, satisfy these conditions. 
This completes the construction of monomial ideals with depth functions of Type II and, therefore, the proof of Theorem \ref{thm.HNTT}.

The following concrete example illustrates the construction of monomial ideals with depth functions of Type II.
 
\begin{example} \label{Type II}
Let $A  = \kk[x,y,z]$ and $I=(x^{d+2},x^{d+1}y,xy^{d+1},y^{d+2},x^dy^2z)$, for some $d \ge 2$.
By \cite[Theorem 4.1]{HH2005} we have
$$\depth A/I^t = \begin{cases} 0 &\text{if $t \le d-1$},\\ 1 &\text{if $t \ge d$}. \end{cases}$$
Let $B = \kk[w,u,v]$. Let $J$ be the integral closure of the ideal
$(w^{3d+3},wu^{3d+1}v,u^{3d+2}v)^3$ or $J = (w^{d+1},wu^{d-1}v,u^dv)$.
By \cite[Example 4.10]{HTT2016} we have
$$\depth B/J^t=\begin{cases} 1 & \text{if $t\le d$},\\ 0 &\text{if $t\ge d+1$}. \end{cases}$$
Let $R = \kk[x,y,z,w,u,v]$. 
By Proposition \ref{additivity}, we have
$$\depth R/((IJ)^t,y-u) = \begin{cases} 1 &\text{if $t\neq d$},\\ 2 &\text{if $t = d$}. \end{cases}$$
Using Proposition \ref{control}, it is easy to check that $z-v$ is a non-zerodivisor modulo $((IJ)^t,y-u)$ for all $t > 0$. Therefore,
$$\depth R/((IJ)^t,y-u,z-v) = \begin{cases} 0 &\text{if $t\neq d$},\\ 1 &\text{if $t = d$}. \end{cases}$$
If we set $S = \kk[x,w,u,v]$ and $Q = (x^{d+2},x^{d+1}u,xu^{d+1},u^{d+2},x^du^2v)J$, which is obtained from $IJ$ by setting $y = u$ and $z = v$, then
$$\depth S/Q^t = \depth R/((IJ)^t,y-u,z-v).$$
Hence, the depth function of $Q$ is of Type II.
\end{example}

Theorem \ref{thm.HNTT} also settles affirmatively a long standing question of Ratliff in \cite[(8.9)]{Rat1983}, that has remained open since 1983.

\begin{question}[Ratliff]
	\label{question.Ratliff}
	Given a finite set $\Gamma$ of positive integer, do there exist a Noetherian ring $S$, an ideal $Q$ and a prime ideal $P \supseteq Q$ in $S$ such that $P$ is an associated prime of $Q^t$ if and only if $t \in \Gamma$?
\end{question}

Specifically, the following corollary is an immediate consequence of Theorem \ref{thm.HNTT}.

\begin{corollary}
	\label{cor.Ratliff}
	Let $\Gamma$ be a set of positive integers which either is finite or contains all sufficiently large integers. There exists a monomial ideal $Q$ in a polynomial ring $S$, with maximal homogeneous ideal $\mm$, such that $\mm \in \Ass(Q^t)$ if and only if $t \in \Gamma$.
\end{corollary}

Corollary \ref{cor.Ratliff}, furthermore, gives a monomial ideal as a counterexample to the following question, which was also due to Ratliff \cite[(8.4)]{Rat1983}. This question was answered nagatively by Huckaba \cite[Example 1.1]{Huck1987}, in which the given ideal was not a monomial ideal.

\begin{question}[Ratliff]
	\label{question.Ratliff2}
	Let $Q$ be an arbitrary ideal in a Noetherian ring $S$. Let $P \supseteq Q$ be a prime ideal such that $P \in \Ass(Q^m)$ for some $m \ge 1$ and $P \in \Ass(Q^t)$ for all $t \gg 0$. Is $P \in \Ass(Q^t)$ for all $t \ge m$?
\end{question}


\section{Symbolic depth functions} \label{sec.symbolic}

Let $Q$ be an ideal in a Noetherian ring $S$. For $t \ge 0$, the $t$-th \emph{symbolic power} of $Q$ is the ideal
$$Q^{(t)} := \bigcap_{\pp \in \Min(Q)} (Q^t_\pp \cap S).$$
In other words, $Q^{(t)}$ is the intersection of the primary components of the minimal associated primes of $Q^t$. We remark here that there is another variant of symbolic powers, in which $\Min(Q)$ is replaced by $\Ass(Q)$, that has also been much investigated. If $Q$ is a radical ideal in a polynomial ring then these definitions agree.
Symbolic powers of homogeneous ideals are much harder to study compared to their ordinary powers. 
This is seen from, for example, the fact that the generators of $Q^{(t)}$ in general cannot be derived merely from 
the generators of $Q$.

Inspired by Theorem \ref{thm.Brod}, one may incline to ask if the symbolic depth function is also necessarily a convergent numerical function; that is, if $\depth S/I^{(t)} = \depth S/I^{(t+1)}$ for all $t \gg 0$. Theorem \ref{HH1} does not apply in this case because the \emph{symbolic Rees algebra} 
$$R_s(Q) := \bigoplus_{t \ge 0}Q^{(t)}$$ 
is not always a standard graded $S$-algebra; it needs not even be finitely generated (see, for instance, \cite{Cut1991, Hun1982,Rob1985}). An application of Theorem \ref{HH1}, when the symbolic Rees algebra of $Q$ is finitely generated, gives us the following result.

\begin{proposition}  \label{Rees}
Let $Q$ be a homogeneous ideal in a polynomial ring $S$. Assume that $R_s(Q)$ is a finitely generated $S$-algebra.
Then, $\depth S/Q^{(t)}$ is an asymptotically periodic function, i.e., it is periodic for $t \gg 0$.
\end{proposition}

By \cite[Theorem 3.2]{HHT2007}, $R_s(Q)$ is finitely generated if $Q$ is a monomial ideal. Therefore, the symbolic depth functions of monomial ideals are asymptotically periodic. For several classes of squarefree monomial ideals, it is known that their symbolic depth functions are actually convergent functions (cf. \cite{CPSTY2015,HKTT2017,KTY2018,TT2012}). 
It was an open question whether the symbolic depth function of any monomial ideal is convergent \cite[p. 308]{HKTT2017}. 

\begin{remark} It is an easy observation that if $\depth S/Q^{(t)} = 0$ for some $t > 0$ then $\depth S/Q^{(t)} = 0$ for all $t > 0$. Therefore, we shall only consider \emph{positive} symbolic depth functions.
\end{remark}

In this section, we survey is a recent result of the second author and H.D. Nguyen \cite{HopT2019a}, which shows that any asymptotically periodic positive numerical function is the symbolic depth function of a homogeneous ideal. In particular, there are plenty monomial ideals whose symbolic depth functions are not necessarily convergent. 

\begin{theorem} \cite[Theorem 6.1]{HopT2019a}
	\label{thm.periodic}
	Let $\kk$ be a field and let $\phi: \NN \rightarrow \NN$ be an asymptotically periodic positive numerical function. Then, there exist a polynomial ring $S$ over a purely transcendental extension of $\kk$ and a homogeneous ideal $Q \subseteq S$ which admits $\phi$ as its symbolic depth function, i.e.,
	$$\depth S/Q^{(t)} = \phi(t) \text{ for all } t \ge 1.$$
\end{theorem}

The proof of Theorem \ref{thm.periodic} is inspired by that of Theorem \ref{thm.HNTT}. The key idea is to construct any asymptotically periodic positive numerical function from basic symbolic depth functions by using closed operations within the class of symbolic depth functions.

Once again, let $A$ and $B$ be polynomial rings over $\kk$ with disjoint sets of variables, and let $R = A \otimes_\kk B$. Let $I \subseteq A$ and $J \subseteq B$ be nonzero proper homogeneous ideals. It follows from Lemma \ref{HoaTam}(i) that 
$$(IJ)^{(t)} = I^{(t)} \cap J^{(t)} = I^{(t)}J^{(t)}.$$
This, together with Lemma \ref{HoaTam}(ii), implies that
$$\depth R/(IJ)^{(t)} = \depth A/I^{(t)}+\depth B/J^{(t)} +1.$$

As in the study of depth function, at this point, we need to find a Bertini-type theorem to get the additivity property of symbolic depth functions. That is, for a given polynomial ring $R = \kk[x_1, \dots, x_n]$ and a homogeneous ideal $K \subseteq R$, we need to find a linear form $f \in R$ such that for all $t \ge 1$, $f$ is a non-zerodivisor of $K^{(t)}$ and if we set $S = R/(f)$ and $Q = (K,f)/(f)$, then
$$S/Q^{(t)} \simeq R/(K^{(t)},f).$$
The first difficulty in finding such a result is that $f$ has to be the same for all symbolic powers $K^{(t)}$, which form an infinite families of ideals. 

The method employed in \cite{HopT2019a} to address this issue is using \emph{generic hyperplane section}. Let $u =\{u_1, \dots, u_n\}$ be a collection of indeterminates and let $R(u) = R \otimes_\kk \kk(u)$, where $\kk(u) = \kk(u_1, \dots, u_n)$ is a purely transcendental extension of $\kk$.
Set  
$$f_u := u_1x_1 + \cdots + u_nx_n.$$
We call $f_u$ a {\em generic linear form}. 
The associated primes of the ideal $(K^{(t)},f_u)$ were already studied in a more general setting in \cite{HKTT2017}. Using results from \cite{HKTT2017}, the following Bertini-type theorem was given in \cite{HopT2019a}.

\begin{proposition} \label{Bertini} \cite[Proposition 5.3]{HopT2019a}
Let $R$ be a polynomial ring over $\kk$ and let $K \subseteq R$ be an ideal with $\depth R/K^{(t)} \ge 2$ for some $t \ge 1$.
Let $S = R(u)/(f_u)$ and $Q = (K,f_u)/(f_u)$.
Then $f_u$ is a regular element on $K^{(t)}R(u)$  and 
$$S/Q^{(t)} \simeq R(u)/(K^{(t)},f_u).$$
\end{proposition}

Proposition \ref{Bertini} has the following consequences on symbolic depth functions.

\begin{corollary} \label{lower}
Let $\phi(t)$ be a symbolic depth function over a field $\kk$ such that $\phi(t) \ge 2$ for all $t \ge 1$.
Then $\phi(t) - 1$ is also a symbolic depth function over a purely transcendental extension of $\kk$.
\end{corollary}

\begin{corollary} \label{additive}
Let $\phi(t)$ and $\psi(t)$ be symbolic depth functions over a field $\kk$.
Then $\phi(t) + \psi(t) - 1$ is a symbolic depth function over a purely transcendental extension of $\kk$.
\end{corollary}

Corollaries \ref{lower} and \ref{additive} particularly show that the operations 
\begin{align*}
\overline \phi(t) & := \phi(t)-1,\\
(\phi*\psi)(t) & := \phi(t)+\psi(t)-1
\end{align*}
are closed in the set of symbolic depth functions with values $\ge 2$ and the set of all symbolic depth functions, respectively.  
 
It is not hard to see that any asymptotically periodic positive numerical function is obtained from finitely many functions of the following types by using the operations $\overline \phi$ with $\phi(t) \ge 2$ for all $t \ge 1$ and $\phi*\psi$:
\begin{enumerate}[Type A:]
\item $1,\dots,1,2,2,\dots$, which is a monotone function converging to 2,
\item $1,\dots,1,2,1,1,\dots$, which has the value 2 at only one position, 
\item $1,1,1,\dots$ or $1,\dots,1,2,1,\dots,1,1,\dots,1,2,1,\dots,1,\dots$, which is a periodic function with a period of the form $1,\dots,1,2,1,\dots,1$, where 2 can be at any position.
\end{enumerate}
The proof of Theorem \ref{thm.periodic} now reduces to showing that all functions of types A, B and C are symbolic depth functions of homogeneous ideals. In fact, any function of types A, B or C is the symbolic depth function of a monomial ideal. This is the most difficult part of the arguments in \cite{HopT2019a}. 

By focusing on monomial ideals, whose symbolic powers are then also monomial ideals, one can invoke a formula of Takayama \cite{Tak2005}, which relates local cohomology modules of a monomial ideal with the reduced homology groups of certain simplicial complexes. 
Since depth can be characterized by the vanishing of the local cohomology modules,
the study of symbolic depth functions can be reduced to the investigation of combinatorial properties of monomial ideals.

To be more precise, let $R = \kk[x_1, \dots, x_n]$ be a polynomial ring over $\kk$ and let $\mm$ be its maximal homogeneous ideal. Let $K \subseteq R$ be a monomial ideal. Note that
$$\depth R/K = \min\{i \mid H_\mm^i(R/K) \neq 0\}.$$

Since $R/K$ has a $\NN^n$-graded structure, 
the local cohomology modules $H_\mm^i(R/K)$ also have a $\ZZ^n$-graded structure.
For $\a \in \ZZ^n$, let $H_\mm^i(R/K)_\a$ denote the degree $\a$ component of $H_\mm^i(R/K)$.
Takayama \cite{Tak2005} gave a formula to relate the dimension and vanishing of $H^i_\mm(R/K)_\a$ to that of the reduced homology groups of certain simplicial complex $\Delta_\a(K)$, which depends on the primary component of $K$. The simplicial complex $\Delta_\a(K)$ is a subcomplex of the Stanley-Reisner simplicial complex $\Delta(K)$ of the squarefree monomial ideal $\sqrt{K}$. Particularly, the facets of $\Delta_\a(K)$ are facets of $\Delta(K)$ if $\a \in \NN^n$. (See \cite{MT2011} for more details and a different interpretation of Takayama's formula.) 

A consequence of Takayama's formula is the following criterion for $\depth R/K \ge 2$.

\begin{proposition} {\rm (cf. \cite[Proposition 1.4]{HopT2019a})}
Let $K$ be an unmixed ideal in $R$. Then $\depth R/K \ge 2$ if and only if $\Delta_\a(K)$ is connected for all $\a \in \NN^n$.
\end{proposition}

On the other hand, if $K$ has a minimal prime $M$ such that $\dim R/M = 2$, then $\depth R/K \le 2$  (see, e.g., \cite[Proposition 1.2.13]{BH1993}). Since $K$ and $K^{(t)}$ share the same minimal primes, in this case, we also have $1 \le \depth R/K^{(t)} \le 2$ for all $t \ge 1$. Hence, the symbolic depth function of $K$ is a 1-2 functions. Note that $\Delta(K) =  \Delta(K^{(t)})$ for all $t \ge 1$.
If we choose $K$ such that $\Delta(K)$ has only one disconnected subcomplex $\Delta'$ whose facets are facets of $\Delta(K)$, then we only need to check when $\Delta_\a(K^{(t)}) = \Delta'$ for all $\a \in \NN^n$ in order to know when $\depth R/K^{(t)} = 1$. 
An instance when this observation applies is given in the following proposition, in which for a monomial ideal with 1-2 symbolic depth function we can test which symbolic powers has depth exactly 2.

\begin{proposition} \label{technical} \cite[Proposition 3.7]{HopT2019a}
\label{lem_key}
Let $R=\kk[x,y,z,u,v]$ be a polynomial ring. Let $M,P,Q$ be primary monomial ideals of $R$ such that 
\begin{align*}
\sqrt{M} & =(x,y,z),\\
\sqrt{P} & =(x,y),\\
\sqrt{Q} & =(z).
\end{align*}
Let $K = M\cap (P,u) \cap (Q,v)$.  Then $\depth R/K^{(t)} \le 2$ and $\depth R/K^{(t)}=2$ if and only if $M^t \subseteq P^t + Q^t$ for all $t > 1$.
\end{proposition}

For the ideal $K$ in Proposition \ref{technical}, we have
$\Delta(K) = \langle \{u,v\},\{z,v\},\{x,y,u\}\rangle,$
which consists of two disjoint facets $\{z,v\}$ and $\{x,y,u\}$ that are connected by the facet $\{u,v\}$. Therefore, $\Delta' = \langle \{z,v\},\{x,y,u\}\rangle$ is the only disconnected subcomplex of $\Delta(K)$ whose facets are facets of $\Delta(K)$. 

Proposition \ref{technical} allows us to construct monomial ideals admitting any given function of Types A, B, C as symbolic depth functions and, thus, completes the proof of Theorem \ref{thm.periodic}. This construction is illustrated in the following examples.

\begin{example} \cite[Lemma 4.2]{HopT2019a} 
	\label{ex.typeA}
	Let $m \ge 2$ be a fixed integer and let $R = \kk[x,y,z,u,v]$. Consider the ideal
	$$K = (x^{2m-2}, y^m, z^{2m})^2 \cap (x^{2m-1}, y^{2m-1}, u) \cap (z,v).$$
	Then, the symbolic depth function of $K$ is of Type A, i.e.,
	$$\depth R/K^{(t)} = \left\{ \begin{array}{ll} 1 & \text{if } t \le m-1, \\ 2 & \text{if } t \ge m.\end{array}\right.$$
\end{example}

\begin{example} \cite[Lemma 4.3]{HopT2019a}
	\label{ex.typeB}
	Let $m \ge 1$ be a fixed integer and let $R = \kk[x,y,z,u,v]$. Consider the ideal
	$$K = (x^{2m}, y^{2m}, xy^{m-1}z, z^{2m})^2 \cap (x^m, y^m, u) \cap (z^{2m+2}, v).$$
	Then, the symbolic depth function of $K$ is of Type B, i.e.,
	$$\depth R/K^{(t)} = \left\{ \begin{array}{ll} 2 & \text{if } t = m, \\ 1 & \text{if } t \not= m.\end{array}\right.$$
\end{example}

For functions of Type C, we first note that the existence of the symbolic depth function $1,1,1,\dots$ is trivial, for example, with $R = \kk[x,y]$ and $K = (x)$. The construction of other symbolic depth functions of Type C is much more subtle because these functions are periodic. (For instance, the construction depends on the period of the given function.) The existence of ideals with symbolic depth functions of Type C is summarized in the following result.

\begin{theorem} \cite[Theorem 4.4]{HopT2019a}
\label{periodic}
Let $m\ge 2$ and $0 \le d < m$ be integers. 
There exists a monomial ideal $K$ in $R=\kk[x,y,z,u,v]$ such that
$$
\depth R/K^{(t)}= \begin{cases}
2 &\text{if $t \equiv d$ modulo $m$},\\
1 &\text{otherwise}.
\end{cases}
$$
\end{theorem}


\section{Open questions} \label{sec.question}

In this section, we discuss open problems and questions related to depth functions and symbolic depth functions that we would like to see answered.

The following question arises naturally from the relationship between depth and projective dimension: 
\begin{quote}
	``\emph{Which numerical functions describe the projective dimension of powers and symbolic powers of homogeneous ideals in polynomial rings?}''
\end{quote}
\noindent This question seems to be very difficult. We could not give the answer even in the following basic situation.

\begin{question}
	\label{quest.pd1}
	Let $Q$ be a homogeneous ideal in a polynomial ring $S$. Suppose that $\pd Q = 1$ and $\pd Q^t = 1$ for all $t \gg 0$. Is it true that $\pd Q^t = 1$ for all $t \ge 1$?
\end{question}

It follows from the Auslander-Buchsbaum formula and Brodmann's result that the projective dimension of powers of an ideal is a convergent function.
Inspired by Theorem \ref{thm.HNTT} and Question \ref{quest.pd1}, we raise the following question.

\begin{question}
	\label{quest.pd}
	Let $g: \NN \rightarrow \ZZ$ be a convergent function such that $g(t) \ge 2$ for all $t \ge 1$. Does there exist a monomial ideal $Q$ in a polynomial ring $S$ such that $g(t) = \pd Q^t$ for all $t \ge 1$?
\end{question}

As a consequence of Theorem \ref{thm.HNTT}, we can give partial answer to Question \ref{quest.pd}.

\begin{corollary}
	\label{cor.PD}
	Let $g: \NN \rightarrow \ZZ_{\ge 0}$ be any convergent numerical function. There exists a monomial ideal $Q$ and a number $c$ such that $\pd Q^t = g(t)+c$ for all $t \ge 1$.
\end{corollary}
\noindent The constant $c$ in Corollary \ref{cor.PD} can be computed as follows. 
Let $m = \max_{t \ge 1} g(t)$. Then $f(t) = m-g(t)$ is a convergent numerical function. Let $n$ be be the number
of variables of a polynomial ring $S$ which contains a homogeneous ideal $Q$ such that 
$\depth S/Q^t = f(t)$ for all $t \ge 1$. Then $\pd Q^t = g(t)+c$ for $c = n-m-1$.
To this end, it is of interest to have an answer to the following question.

\begin{question} 
What is the smallest number of variables of a
polynomial ring which contains a homogeneous ideal with a given depth function $f(t)$?
\end{question}
\noindent Note that the proof of Theorem \ref{thm.HNTT} uses a large number of variables compared to the values of $f(t)$. 

In making use of Corollaries \ref{lower} and \ref{additive}, the ideals constructed in Theorem \ref{thm.periodic} are non-monomial ideals in polynomial rings over purely transcendental extensions of the given field $\kk$. 
Using the theory of specialization \cite{Krull1948,Nhi2007,NhiTrung1999}, we can construct such ideals in polynomial rings over any uncountable field. This is because the Bertini-type result, Proposition \ref{Bertini}, holds without having to go to purely transcendental extensions of the ground field; see \cite[Proposition 5.8]{HopT2019a}.
This leads us to the following question.

\begin{question} 
Given a field $\kk$ and an asymptotically periodic positive numerical function $\phi(t)$, 
do there exist a polynomial ring $S$ over $\kk$ and a {\em monomial ideal} $Q \subset S$ such that $\depth S/Q^{(t)} = \phi(t)$ for all $t \ge 1$?
\end{question}
\noindent The analogous question for the depth function of homogeneous ideals has a positive answer by Theorem \ref{thm.HNTT}. 

Again, due to the relationship between depth and projective dimension and inspired by Theorem \ref{thm.periodic}, we raise the following question on projective dimension of symbolic powers.

\begin{question}
	\label{quest.pds}
	Let $g: \NN \rightarrow \NN$ be an asymptotically periodic function such that $g(t) \ge 2$ for all $t \ge 1$. Does there exist a monomial ideal $Q$ in a polynomial ring $S$ such that $g(t) = \pd Q^{(t)}$ for all $t \ge 1$?
\end{question}

As a consequence of Theorem \ref{thm.periodic}, we obtain a partial answer to Question \ref{quest.pds}.

\begin{corollary}
	\label{cor.pds}
	Let $\phi: \NN \rightarrow \ZZ_{\ge 0}$ be an asymptotically periodic numerical function. Let $\kk$ be a field and let $m = \max_{t \ge 1} \phi(t)$. Then, there exist a positive integer $c$, a polynomial ring $S$ in $m+c+2$ variables over a purely transcendental extension of $\kk$, and a homogeneous ideal $Q$ in $S$ such that
	$\pd Q^{(t)} = \phi(t) + c \text{ for all } t \ge 1.$
\end{corollary}

Similarly to Corollary \ref{cor.PD}, the constant $c$ in Corollary \ref{cor.pds} is determined by the number of variables of a polynomial ring $S$ which contains a homogeneous ideal $Q$ with the given symbolic depth function.

The proof of Theorem \ref{thm.periodic} uses a large number of variables.  However, all constructed examples of symbolic depth functions of types A, B, C (except $1,1,1,\dots$) are ideals of height 2 in polynomial rings in 5 variables.
It is naturally of interest to consider the following question.

\begin{question} \label{variables}
Let $\phi(t)$ be an asymptotically periodic positive numerical function and $m = \max_{t\ge1} \phi(t)$.
Does there exist a polynomial ring $S$ in $m+3$ variables that contains a height 2 homogeneous ideal $Q$
such that $\depth S/Q^{(t)} = \phi(t)$ for all $t \ge 1$?
\end{question}

Theorem \ref{thm.periodic} classifies a large class of symbolic depth functions. It remains an open problem to determine if Theorem \ref{thm.periodic} indeed covers all symbolic depth functions. 

\begin{question} 
Does there exist a homogeneous ideal whose symbolic depth function is not asymptotically periodic?
\end{question}
\noindent According to Proposition \ref{Rees}, if such an ideal existed, its symbolic Rees algebra would have to be non-Noetherian.
To find non-Noetherian symbolic Rees algebras is a difficult problem that is related to Hilbert's fourteenth
problem; see, for instance, \cite{Rob1985}. 
To the best of our knowledge, there are only examples of non-Noetherian symbolic Rees algebras for one-dimensional ideals  (cf. \cite{Cut1991, Hun1982, Rob1985}). In this case, we have $\depth S/I^{(t)} = 1$ for all $t \ge 1$, whence the symbolic depth function is a constant function.

It was shown in \cite{HopT2019a,MN2018} that the symbolic depth function of a squarefree monomial ideal $Q$ is \emph{almost} non-increasing, in the sense that $\depth S/Q^{(s)} \le \depth S/Q^{(t)}$ for $s \gg t$. There are examples of ideals generated by squarefree monomials of degrees $\ge 3$ whose symbolic depth functions need not be monotone \cite{HopT2019a}.

\begin{question}
	\label{quest.monotone}
	Is the symbolic depth function of the edge ideal of a graph a non-increasing function?
\end{question}
\noindent The analogous question for the depth function of the edge ideal of a graph is also an open question (cf. \cite{HH2005, HQ2015}). Note that the depth function of a squarefree monomial ideal in general needs not be non-increasing; see \cite{HS2015, KSS2014}.

Beside powers and symbolic powers of an ideal, the integral closures of powers have been extensively investigated. It is also a classical result of Brodmann \cite{Brod1979} that for an ideal $Q$ in a Noetherian ring $S$, $\depth S/\overline{Q^t}$ is asymptotically a constant function, i.e., the function $\depth S/\overline{Q^t}$ is a convergent numerical function.

\begin{question}
	\label{quest.int}
	For which convergent numerical function $f: \NN \rightarrow \ZZ_{\ge 0}$ does there exist a homogeneous ideal $Q$ in a polynomial ring $S$ such that $\depth S/\overline{Q^t} = f(t)$ for all $t \ge 1$?
\end{question}

\noindent The monomial generators of $\overline{Q^t}$ can be derived from that of $Q$ by combinatorial means; see, for instance, \cite{HT2019}. This fact was used in \cite{HopT2019a} to examine the depth of integrally closed symbolic powers of monomial ideals.



\begin{thebibliography}{999}
\bibitem{AB1956} M. Auslander and D. A. Buchsbaum. Homological dimension in Noetherian rings. Proc. Nat. Acad. Sci. U.S.A., 42:36--38, 1956.

\bibitem{BHH2014} S. Bandari, J. Herzog, and T. Hibi. Monomial ideals whose depth function has any given number of
strict local maxima. Ark. Mat., 52(1):11--19, 2014.

\bibitem{Brod1979} M. Brodmann. The asymptotic nature of the analytic spread. Math. Proc. Cambridge Philos. Soc.,
86(1):35--39, 1979.

\bibitem{BH1993} W. Bruns and J. Herzog. \emph{Cohen-Macaulay rings}, volume 39 of \emph{Cambridge Studies in Advanced Mathematics}. Cambridge University Press, Cambridge, 1993.

\bibitem{CPSTY2015} Constantinescu, M. R. Pournaki, S. A. Seyed Fakhari, N. Terai, and S. Yassemi. Cohen-Macaulayness
and limit behavior of depth for powers of cover ideals. Comm. Algebra, 43(1):143--157, 2015.

\bibitem{Cut1991} S. D. Cutkosky. Symbolic algebras of monomial primes. J. Reine Angew. Math., 416:71--89, 1991.

\bibitem{HNTT2020PAMS} H. T. H\`a, H. D. Nguyen, N. V. Trung, and T. N. Trung. Depth functions of powers of homogeneous
ideals. Proc. Amer. Math. Soc., 149: 1837--1844, 2021.

\bibitem{HS2015} H. T. H\`a and M. Sun. Squarefree monomial ideals that fail the persistence property and non-increasing
depth. Acta Math. Vietnam., 40(1):125--137, 2015.

\bibitem{HT2019} H. T. H\`a and N. V. Trung. Membership criteria and containments of powers of monomial ideals. Acta
Math. Vietnam., 44(1):117--139, 2019.

\bibitem{HTT2016} H. T. H\`a, N. V. Trung, and T. N. Trung. Depth and regularity of powers of sums of ideals. Math. Z.,
282(3-4):819--838, 2016.

\bibitem{HQ2015} J. Herzog and A. Asloob Qureshi. Persistence and stability properties of powers of ideals. J. Pure Appl.
Algebra, 219(3):530--542, 2015.

\bibitem{HH2005} J. Herzog and T. Hibi. The depth of powers of an ideal. J. Algebra, 291(2):534--550, 2005.

\bibitem{HHT2007} J. Herzog, T. Hibi, and N. V. Trung. Symbolic powers of monomial ideals and vertex cover algebras.
Adv. Math., 210(1):304--322, 2007.

\bibitem{HKTT2017} L. T. Hoa, K. Kimura, N. Terai, and T. N. Trung. Stability of depths of symbolic powers of Stanley-Reisner ideals. J. Algebra, 473:307--323, 2017.

\bibitem{HT2010} L. T. Hoa and N. D. Tam. On some invariants of a mixed product of ideals. Arch. Math. (Basel), 94(4):327--337, 2010.

\bibitem{Huck1987} S. Huckaba. On linear equivalence of the P-adic and P-symbolic topologies. J. Pure Appl. Algebra,
46(2-3):179--185, 1987.

\bibitem{Hun1982} C. Huneke. On the finite generation of symbolic blow-ups. Math. Z., 179(4):465--472, 1982.

\bibitem{KSS2014} T. Kaiser, M. Stehl\'ik, and R. \v{S}krekovski. Replication in critical graphs and the persistence of monomial ideals. J. Combin. Theory Ser. A, 123:239--251, 2014.

\bibitem{KTY2018} K. Kimura, N. Terai, and S. Yassemi. The projective dimension of the edge ideal of a very well-covered
graph. Nagoya Math. J., 230:160--179, 2018.

\bibitem{Krull1948} W. Krull. Parameterspezialisierung in Polynomringen. Arch. Math., 1:56--64, 1948.

\bibitem{MST2018} K. Matsuda, T. Suzuki, and A. Tsuchiya. Nonincreasing depth functions of monomial ideals. Glasg.
Math. J., 60(2):505--511, 2018.

\bibitem{MT2011} N. C. Minh and N. V. Trung. Cohen-Macaulayness of monomial ideals and symbolic powers of Stanley-
Reisner ideals. Adv. Math., 226(2):1285--1306, 2011.

\bibitem{MN2018} J. Monta\~no and L. N\'u\~nez Betancourt. Splittings and symbolic powers of squarefree monomial ideals.
To appear in Int. Math. Res. Notices.

\bibitem{HopT2019a} H. D. Nguyen and N. V. Trung. Depth functions of symbolic powers of homogeneous ideals. Invent.
Math., 218(3):779--827, 2019.

\bibitem{HopT2019b} H. D. Nguyen and N. V. Trung. Correction to: Depth functions of symbolic powers of homogeneous
ideals. Invent. Math., 218(3):829--831, 2019.

\bibitem{Nhi2007} D. V. Nhi. Specializations of direct limits and of local cohomology modules. Proc. Edinb. Math. Soc.
(2), 50(2):459--475, 2007.

\bibitem{NhiTrung1999} D. V. Nhi and N. V. Trung. Specialization of modules. Comm. Algebra, 27(6):2959--2978, 1999.

\bibitem{Rat1983} L. J. Ratliff, Jr. A brief survey and history of asymptotic prime divisors. Rocky Mountain J. Math.,
13(3):437--459, 1983.

\bibitem{Rob1985} P. C. Roberts. A prime ideal in a polynomial ring whose symbolic blow-up is not Noetherian. Proc.
Amer. Math. Soc., 94(4):589--592, 1985.

\bibitem{Serre1965} J.-P. Serre. \emph{Alg\`ebre locale. Multiplicit\'es}, volume 11 of \emph{Cours au Coll\`ege de France, 1957--1958, r\'edig\'e par Pierre Gabriel. Seconde \'edition, 1965. Lecture Notes in Mathematics}. Springer-Verlag, Berlin-New York, 1965.

\bibitem{Tak2005} Y. Takayama. Combinatorial characterizations of generalized Cohen-Macaulay monomial ideals. Bull.
Math. Soc. Sci. Math. Roumanie (N.S.), 48(96)(3):327--344, 2005.

\bibitem{TT2012} N. Terai and N. V. Trung. Cohen-Macaulayness of large powers of Stanley-Reisner ideals. Adv. Math.,
229(2):711--730, 2012.

\end{thebibliography}
\end{document}